\documentclass[12pt, leqno]{article}
\usepackage{amsmath}
\usepackage{amsfonts}
\usepackage{amssymb}
\usepackage{graphicx}
\usepackage[monochrome]{color}
\newtheorem{thm}{Theorem}[section]
\newtheorem{cor}[thm]{Corollary}
\newtheorem{lem}[thm]{Lemma}

\numberwithin{equation}{section}

\newcommand{\RR}{\mathbb{R}}
\newcommand{\ren}{\mathbb{R}^n}

\newcommand{\wu}{\widetilde{u}}
\newcommand{\ve}{\varepsilon}

\newcommand{\dint}{\displaystyle\int}
\newcommand{\diint}{\displaystyle\iint}
\renewcommand{\dfrac}{\displaystyle\frac}
\def\qed{\,\unskip\kern 6pt \penalty 500
\raise -2pt\hbox{\vrule \vbox to8pt{\hrule width 6pt
\vfill\hrule}\vrule}\par}
\begin{document}
\title{\textbf{Nonlinear porous medium flow  \\
with fractional potential pressure}}
\author{\Large Luis Caffarelli \footnote{caffarel@math.utexas.edu}
~and~ Juan Luis Vazquez\footnote{juanluis.vazquez@uam.es}\\}
\date{ }
\maketitle

\begin{abstract}
We study a porous medium equation with nonlocal diffusion effects given by an inverse fractional Laplacian operator:
$$
\partial_t u-\nabla\cdot(u\nabla p)=0, \quad p=(-\Delta)^{-s}u,\qquad 0<s<1.
$$
We pose the problem for $x\in \ren$ and $t>0$ with bounded and compactly supported initial data, and prove existence of weak and bounded solutions that propagate with finite speed, a property that is not shared by other fractional diffusion models.

\end{abstract}


\section{Introduction}
\label{sec.intro}

We study a nonlinear diffusion model with nonlocal effects  described by the system
\begin{equation}\label{eq1}
\partial_t u=\nabla\cdot(u\nabla p), \quad p={\cal K}(u).
\end{equation}
Here, $u$ is a function of the variables $(x,t)$ to be thought of
as a density or concentration, and therefore nonnegative, while
$p$ is the pressure, which is related to $u$ via a linear \color{blue} positive \normalcolor operator
${\cal K}$, which we assume to be the inverse of a fractional Laplacian.
To be specific, the problem is posed for $x\in \RR^n$, $n\ge 1$, and $t>0$, and we
give initial conditions
\begin{equation}\label{eq.ic}
u(x,0)=u_0(x), \quad x\in \RR^n,
\end{equation}
where $u_0$ is a nonnegative and bounded function in $\RR^n$ with compact support
or fast decay at infinity.

The model arises from the consideration of a continuum, say, a
fluid or a population, represented by a  density distribution $u(x,t)$ that
evolves with time following a velocity field ${\bf v}(x,t)$,
according to the continuity equation
$$
u_t+\nabla\cdot(u\, {\bf v})=0.
$$
We now assume that $\bf v$ derives from a potential, ${\bf
v}=-\nabla p$, as happens for instance in fluids in porous media
according to Darcy's law, and in that case $p$ denotes the pressure. But
potential velocity fields are also found in many other instances, like
Hele-Shaw cells.

We still need a closure relation to relate $u$ and $p$. In the
case of gases in porous media, as modeled in the 1930's by Leibenzon and Muskat
\cite{Leib, Mu37}, the closure relation takes the form of a state law: $p=f(u)$, where
$f$ is a nondecreasing scalar function, which is linear when the
flow is isothermal, and a higher power of $u$ if it is adiabatic.
The linear relationship happens also in the
simplified description of water infiltration in an almost
horizontal soil layer according to Boussinesq. See \cite{Vapme}
for a description of these and other applications. Summing up, we get the standard porous
medium equation, $u_t=c\Delta (u^2)$, or more generally, $u_t=\Delta (u^m)$ with $m>1$.

In this paper we propose to consider  the case where $p$ es
related to $u$ through a  linear fractional potential operator, ${\cal
K}=(-\Delta)^{-s}$ with kernel $K(x,y)= c|x-y|^{-(n-2s)}$ (i.\,e., a
Riesz operator, cf. \cite{Land, Stein}, see also Appendix for precise definitions and some comments).
The interest in using fractional Laplacians in modeling diffusive processes has a wide literature,
 especially when one wants to model long-range diffusive interaction,
and this interest has been activated by the recent progress in the
mathematical theory. This literature is mostly elliptic, cf. \cite{ACS, CSS, Silv}, but cf. works like \cite{AC09, P4} for related parabolic problems.

More generally, it could be assumed that ${\cal K}$ is an operator of
integral type defined by convolution on all of $\RR^n$, with the
assumptions that is positive and symmetric. The fact the ${\cal K}$ is a
homogeneous operator of degree $2s$, $0<s<1$, will be important in
the proofs given below.  An interesting variant would be ${\cal K}=(-\Delta+
cI)^{-s}$. We are not exploring here such extensions.

\medskip

\noindent {\sc  Extreme cases.} (1) If we take in our model $s=0$, so that ${\cal K}=$ the
identity operator, we get the standard porous medium equation with $m=2$,
whose behavior is well-known, see  \cite{Ar,Vapme} for the mathematical
theory and the applications.

\noindent (2) In the other end of the $s$ interval, when $s=1$ and
we take ${\cal K}=(-\Delta)^{-1}$, we get
\begin{equation}
u_t=\nabla u\cdot \nabla p - u^2, \quad -\Delta p=u.
\end{equation}
In one dimension this leads to $u_t= u_x p_x - u^2, p_{xx}=-u.$ In
terms of $v=-p_x=\int u\,dx$ we have
$$
v_t=up_x+ c(t)=-v_xv+c(t),
$$
For $c=0$ this is the Burgers equation $v_t+vv_x=0$ which generates shocks
in finite time if we allow for $u$ to have two signs.

As a related precedent we may mention the model studied by Lions and Mas-Gallic \cite{LMG},
who are interested in the regularization of the velocity field in the standard porous
medium equation by means of a convolution kernel. They  get a system  which is formally like ours, but there is a big difference in the study since they assume the kernel to be smooth and integrable, and in fact an approximation of the Dirac delta, in other words, a short-range interaction. Since the kernel of the fractional
operator \color{blue} $(-\Delta)^{-s}$ \normalcolor is $k(x,y)\sim \,|x-y|^{-(n-2s)}$, i.\,e., a long-range interaction, we are far away from that situation, but it  may serve as a previous regularization step.

A model from superconductivity  arises in recent work  by Ambrosio and Serfaty \cite{AmSr} describing the evolution of the vortex-density in superconductor modeling. The system is similar to our system with ${\cal
K}=(-\Delta)^{-s}$, $s=1$ and their mathematical tools are quite different. \color{blue} On the other hand, the equation with $s=1/2$ has been proposed by Head \cite{H} as the equation of motion of a dislocation continuum, and then $u$ is the dislocation density and the space dimension is $n=1$. The mathematical investigation of this case is performed by Biler et al. in \cite{BKM}, though in terms of the integrated equation $v_t+|v_x|\Lambda(v)=0$, where $\Lambda$ is the L\'evy operator of order 1, which is equivalent to $(-\partial^2_{xx})^{1/2}$. \normalcolor

Models of this kind arise in other contexts.  Some variants will be indicated at the end of the paper.

\medskip

\noindent {\bf Organization of the paper.} Section \ref{sec.estimates1} derives the basic estimates in a formal way. The proof of existence of a weak solution proceeds by approximation, whereby the degeneracy of the equation is eliminated, diffusion is added and the kernels are regularized. This is technically delicate, so we first prove existence of weak solutions for the approximate problems posed in bounded domains in  Section \ref{sec.ex1}, and at the time basic estimates are rigorously derived. In Section \ref{sec.ex2} weak solutions of the original problem are constructed in the whole space by  passage to the limit  after a tail control step based on a novel  argument with so-called suitable ``true upper barriers'', cf. Theorem \ref{thm.ex}. Such barrier method in new and turns out to be well adapted to obtain comparison results in the presence of nonlocal operators.

We then establish the main properties of the solutions: Section \ref{sec.fp} establishes the property of finite propagation, which is a main feature  of porous media equations and gives rise to the appearance of a free boundary. We discuss the persistence of positivity in Section \ref{sec.pers}.

Two sections close the paper: the Appendix,  Section  \ref{appendix}, gathers some useful definitions. A final Section \ref{sec-comm} contains comments on variants,  extensions or ongoing work on the topic of this paper: this refers in particular to the pending questions of uniqueness, smoothness or asymptotic behaviour.

\medskip

\noindent {\bf Notation.}  We will use the notation $L_s=(-\Delta)^{s}$ with $0<s<1$ for the
fractional powers of the Laplace operator defined on smooth functions in $\RR^n$ by Fourier transform and extended in a natural way to functions in the Sobolev space $H^{2s}(\RR^n)$.  Technical reasons imply that in one space dimension the restriction $s<1/2$ will be observed. The inverse operator is denoted by ${\cal K}_s=(-\Delta)^{-s}$ and can be realized by convolution
$$
{\cal K}_s=K_s\star u, \qquad K_s(x)=c(n,s)|x|^{-(n-2s)}.
$$
as described in the appendix. ${\cal K}_s$ is a positive self-adjoint operator. We will write ${\cal H}_s={\cal K}_s^{1/2}$ which has kernel $K_{s/2}$. The subscript $s$ will be omitted when $s$ is fixed and known. For functions that depend on $x$ and $t$, convolution is applied for every fixed $t$ with respect to the space variables. We then use the abbreviated notation $u(t)=u(\cdot,t)$.

\section{Basic estimates}
\label{sec.estimates1}

The existence theory of weak solutions needs a lengthy process based on several approximations and passage to the limit that may obscure to the reader the main properties of the solutions. These are however very clear from usual considerations in mathematical physics, and the authors think that it useful for the reader to have them in mind as a goal. Therefore, we do at this stage formal calculations, assuming that $u\ge 0$  satisfies the required smoothness and integrability assumptions and decreases fast enough as $|x|\to \infty$. The calculations are to be justified later by the approximation process. We fix $s\in (0,1)$ and put ${\cal K}=(-\Delta)^{-s}$ and  ${\cal H}={\cal K}^{1/2}$.

 \noindent $\bullet$
Conservation of mass
\begin{equation}
\frac{d}{dt}\int_{\ren} u(x,t)\,dx=\int_{\ren} \nabla \cdot (u\,{\cal K}u)\,dx=0.
\end{equation}

\medskip

\noindent $\bullet$ First energy estimate:
\begin{equation}
\frac{d}{dt}\int_{\ren} u(x,t)\log u(x,t)\,dx=-\int_{\ren} (\nabla u\cdot \nabla {\cal K}u)\,dx=-\int_{\ren} |\nabla {\cal H}u|^2\,dx\,,
\end{equation}
where we use the fact that ${\cal K}={\cal H}^2$, and ${\cal H}$ is a positive self-adjoint operator that commutes with the gradient.

\medskip

\noindent $\bullet$ Second energy estimate:
\begin{equation}
\begin{array}{l}
\dfrac12 \,\dfrac{d}{dt}\int_{\ren} |{\cal H}u(x,t)|^2\,dx=\int_{\ren} {\cal H}u\,({\cal H} u)_t\,dx=\\
\dint_{\ren} {\cal K}u\,u_t\,dx=\dint_{\ren} ({\cal K}u)\,\nabla\cdot(u\nabla {\cal K}u)\,dx=-\dint_{\ren} u|\nabla {\cal K}u|^2\,dx.
\end{array}
\end{equation}

\noindent $\bullet$ $L^\infty $ estimate. We prove that the
$L^\infty$ norm does not increase in time.

\noindent {\sl Sketch of the proof.} At a point of maximum of $u$ at time
$t=t_0$, say $x=0$,  we have
$$
u_t=\nabla u\cdot \nabla P + u\,\Delta {\cal K}(u).
$$
The first term is zero, and for the second we have $-\Delta {\cal K}= L$
where $L=(-\Delta)^q$ with $q=1-s$ so that
$$
\Delta {\cal K}u(0)=-Lu(0)=-c\int_{\ren} \frac{u(0)-u(y)}{|y|^{n+2(1-s)}}\,dy \le
0\,,
$$
where $c(s,n)>0$.

\medskip

\noindent $\bullet$ Conservation of positivity: $u_0\ge 0$ implies that $u(t)\ge 0$
for all times. The argument is similar.

\medskip

\noindent $\bullet$ We  derive next the $L^p$ estimates, $1<p<\infty$ :
\begin{eqnarray*}
&\dfrac{d}{dt}\dint u^p(x,t)\,dx=p\int_{\ren} u^{p-1}u_t\,dx=-p(p-1)\int_{\ren} u^{p-1}\nabla u \cdot \nabla {\cal K}u\,dx =\\
&(p-1)\diint (u(x)^{p}-u(y)^{p})\,\Delta K_s(x-y) u(y)\,dxdy=\\
&-(p-1)\diint (u(x)^{p}-u(y)^{p})\,\Delta K_s(x-y) (u(x)-u(y))\,dxdy\\
&+ (p-1)\diint (u(x)^{p}-u(y)^{p})\,\Delta K_s(x-y) u(x)\,dxdy=-I_1+ I_2.
\end{eqnarray*}
Since $u\ge0$  and  $\Delta K_s\ge 0$ for $0<s<1$,  the first integral $I_1$ is positive. But symmetry means that
$I_2=I_1/2$. We conclude that $\int u^p\,dx$ is decreasing in time. See a related calculation in  \cite{CVss} and \cite{CoCo03}.

\medskip

\ $\bullet$ A standard comparison result for parabolic equations does not seem to work. This is one of the main technical difficulties in the study of this equation. In fact, we will find special situations where some comparison holds. Some partial comparison allows us to prove two main results of the paper.

\section{Existence I. Smooth approximate solutions}
\label{sec.ex1}

We want to  solve the initial-value problem for the equation
\begin{equation}\label{eq}
\partial_t u=\nabla (u\,\nabla {\cal K}u), \qquad {\cal K}=(-\Delta)^{-s},
\end{equation}
posed in $Q=\ren\times (0,\infty)$  or at least $Q_T=\ren\times (0,T) $, with parameter $ 0<s<1$. We will take  initial data $u_0(x)\ge 0$, $u_0\in L^1(\RR^n).$ We  assume mostly for technical convenience that $u_0$ is bounded, and in the next section we will also impose decay conditions as $|x|\to\infty$.

We want to obtain a suitable weak solution $u(x,t)$ defined in $Q$. We approach this problem by a process that consists of regularization, elimination of the degeneracy, and reduction of the space domain. Once the approximate problems are solved, estimates are obtained that allow to pass to the limit step by step in all the approximations to obtain in the end a weak solution of the original problem.

\noindent {\bf Definition.} We say that $u$ is a weak solution of equation (\ref{eq})
in $Q_T=\RR^n\times (0,T)$ with initial data $u_0\in L^1(\RR^n)$ if $u\in L^1(Q_T)$, ${\cal K}(u)\in L^1(0,  T: W^{1,1}_{loc}(\RR^n))$, and $u\,\nabla{\cal K}(u)\in L^1(Q_T)$,
and the identity
\begin{equation}
\iint u\,(\phi_t-\nabla {\cal K}(u)\cdot\nabla\phi)\,dxdt+ \int
u_0(x)\,\phi(x,0)\,dx=0
\end{equation}
holds for all continuous test
functions $\phi$ in $Q_T$ such that $\nabla_x\phi$ is continuous, and
$\phi$ has compact support in the space variable and vanishes near
$t=T$.

\medskip

\noindent \noindent {\bf \ref{sec.ex1}.1.} The modifications that we use as  a starting point are as follows:
regularization is done by adding Laplacian term plus a kernel smoothing; the degeneracy is eliminated by raising  the $u=0$ level in the diffusion coefficient. Specifically, we  take small numbers $\ve,\delta,\mu\in (0,1)$ and consider the equation $(E(\ve,\delta,\mu))$:
\begin{equation}
u_t=\delta \Delta u + \nabla\cdot(d(u)\nabla {\cal K}_\ve(u)),
\end{equation}
posed in $Q_{T,R}=\{x\in B_R(0), 0<t<T\}$. A simple  option for $d(u)$ is  ${d}(u)=u+\mu$ with a small $\mu>0$.
 Another option would be $d(u)=\mu$ for $0\le u\le \mu$ and $d(u)=u$ for $u\ge \mu$ (we prefer the former one). Besides, the equation is posed in the spatial domain $B_R=B_R(0)$ for $0<t<T$. We also take initial conditions
\begin{equation}
u(x,0)=\hat u_0(x) \qquad x\in B_R(0),
\end{equation}
where  $\hat u_0=u_{0,\ve,R}$ is a nonnegative, smooth and bounded approximation
of the initial data $u_0\ge 0$. Finally, we take boundary data
\begin{equation}
u(x,t)=0 \qquad \mbox{for } \ |x|=R, \ t\ge 0.
\end{equation}
Let us explain now a convenient approximation of the kernel:  ${\cal K}_\ve u(t)=\zeta_\ve\star u(t)$, and $\zeta_\ve$ is a smooth approximation  of the Riesz kernel $k_s(x)=c|x|^{-(n-2s)}$ corresponding to the inversion of the $s$-Laplacian on $\RR^n$. In our implementation it acts on the extension of $u(x,t)$ to the whole domain $x\in \RR^n$, which is done in the natural way, i.\,e., putting $u=0$ for $|x|\ge R$ and $t>0$.
This approximation process is a bit similar to the approximation of the porous media performed
by Lions and Mas-Gallic in \cite{LMG}, but their kernel was just a mollifying kernel
representing short-range effects and the consequences of a long-range kernel with slow decay at infinity are quite different.

The existence and uniqueness of a solution $u(x,t) =u_{\ve,\delta,\mu,R}(x,t)$ for the model  is then  more or less standard, and the solution is smooth. In the weak formulation we have
\begin{equation}
\iint u\,(\phi_t-\delta \Delta \phi)\,dxdt -\iint d(u)(\nabla {\cal K}_\ve(u)\cdot\nabla\phi)\,dxdt+
\int u_0(x)\,\phi(x,0)\,dx=0
\end{equation}
with double integrals in $Q_{T,R}=B_R(0)\times (0,T)$, and valid
for every smooth $\phi$ that vanishes at the lateral boundary and for all large $t$.
It is also clear that a priori estimates like the ones
in the previous section apply to this model. In particular we have

\noindent (E1.a) An easy computation gives the decay of the ''total mass''
\begin{equation}
\int_{B_R} u(x,t)\,dx\le \int_{B_R} \hat u_0(x)\,dx.
\end{equation}
Of course, we lose the expected mass conservation of the original problems because of the zero Dirichlet conditions of our present problem, but the estimate is still useful.

\medskip

\noindent (E1.b) We also need conservation of nonnegativity, which is easy.

\medskip

\noindent (E1.c) The $L_x^\infty$ bound is conserved, $0\le u(x,t)\le \|\hat u_0\|_\infty$, and the argument is as in the previous section. As a consequence, the solutions $u(\cdot,t)$ at time $t$ belong to all $L^p(B_R)$ spaces with norm that is independent of the parameters $\delta, \ve, \mu$ and $R$.

\medskip

\noindent (E1.d)  We now introduce a version of the  first energy inequality of previous section. We
select the function $F$ defined by the conditions $F''(u)=1/d(u)$ and $F(0)=F'(0)=0$.
After some integrations by parts we will  get
\begin{equation}
\dfrac{d}{dt}\int_{B_R} F(u)\,dx=-\delta \int_{B_R} \frac{|\nabla u|^2}{d(u)}\,dx -\int_{B_R} |\nabla {\cal  H}_\ve u|^2\,dx\,,
\end{equation}
where ${\cal  H}_\ve={\cal  K}_\ve^{1/2}$ and we have used the fact that $F'(u)=0$ on $\Sigma=\{|x|=R\}\times [0,T]$ to annihilate the boundary term in the integration by parts. This formula implies that for all $0<t<T$:
\begin{equation}
\int_{B_R} F (u(t))\,dx+\delta \int_0^t\int_{B_R} \frac{|\nabla u|^2}{d(u)}\,dx dt+\int_0^t\int_{B_R} |\nabla {\cal  H}_\ve u|^2\,dxdt=\int F(\hat u_0)\,dx.
\end{equation}
This implies estimates for $|\nabla {\cal  H}_\ve u|^2$ \ and \ $\delta |\nabla u|^2/d(u)$ \ in \ $ L^1_{x,t}(Q_{T,R})$ and the bounds for such norms are independent of $\ve$, $\delta$, $R$, and $T$. They  do depend on $\mu>0$ through
the value of $F=F_\mu$. Indeed, the explicit formula for $F(u) $ is defined for $u>0$ as:
\begin{equation}
F_\mu(u)=(u+\mu)\log(1+(u/\mu))-u, \quad  F_\mu'(u)=\log(1+(u/\mu)).
\end{equation}
This offers difficulties when $\mu\to 0$. (In case we take the second option, we would have
$$
F_\mu(u)=(1/2\mu)u^2 \quad \mbox{for }  \ 0\le u\le \mu, \qquad F_\mu(u)=u\log(u/\mu)+(\mu/2) \quad \mbox{otherwise},
$$
which is no better)

\noindent {\bf Note.} We take as ${\cal K}_\ve$ the operator obtained by convolution with a standard mollification of $k_s$ of the form $\zeta_\ve=k_s\star \rho_\ve$ where $\rho_\ve(x)=\ve^{-n}\rho(x/\ve)$ and $\rho$ is a $C^\infty_c(\RR^n)$, $\rho\ge0$, $\rho$ radially symmetric and decreasing; moreover, if $\rho=\sigma\star\sigma$ where $\sigma$ has the same properties. Then, we can write ${\cal H}_\ve$ as the operator with kernel $k_{s/2}\star \sigma_\ve$.
\medskip

\noindent {\bf \ref{sec.ex1}.2.}  We have to pass to the limit in four parameters: $\delta, \ve, \mu$ and $R$. The last two limits are the most delicate. In order to examine the convergence arguments, we can  try to pass to the limit $\ve\to 0$ as a next step to obtain a solution of the equation
\begin{equation}\label{appr.eq}
u_t= \delta\,\Delta u+ \nabla\cdot(d(u)\nabla {\cal K}(u)),
\end{equation}
with same  initial and boundary data.  Using in a  precise way the above estimates we get convergence of $u_{\ve}\to u$ in $L^\infty_t(L^1_x)$ weak. This is enough for the limit of
the first integral of the formula of weak solution. The second integral contains the product $u\,\nabla {\cal K}u$ and we need to study better the consequences of the estimates:

(i)  Since ${\cal K}(u)={\cal H}({\cal H}(u))$ and $\nabla {\cal K}(u)=\nabla {\cal H}({\cal H}(u))={\cal H}(\nabla ({\cal H} u))$, we derive from the bound for $\nabla {\cal H}u$ in $L^2_{x,t}$ the needed estimates for ${\cal K}(u)$. Recall that ${\cal H}(u)$ and ${\cal K}(u)$ are defined on all of $\RR^n$. We recall that $\nabla {\cal H}(u)$ is a "derivative of order $1-s$ of $u$", and since $u$ is bounded, $u\in L^\infty_{x,t}$, we conclude that $u\in L^2_t H^{1-s}_{x,loc}$. By potential theory, it is then clear that ${\cal K}(u)\in L^2_tH_x^{1+s}$.

(ii) We also need  some continuity in time. Using the equation, which expresses $u_t$ as the divergence of $u\,\nabla {\cal K}u,$  together with the boundedness of $u$ and the bound for ${\cal K}(u)$ in $H_x^{1+s}$, we conclude that $u_t\in L_t^2 (H_x^{-1+s})$. Hence, the family  of approximations to $u$ is relatively compact in the sense of the parabolic compactness theorems by Aubin and Simon, cf. \cite{Aubin}, \cite{Simon}. This means that a limit $u$ exists (along suitable subsequences) and $u\in C([0,T]:L^2(B_R))$.

All of this is used in passing to the limit of the term $\iint u (\nabla {\cal K}u)\nabla \phi\,dxdt$ as follows: we have the convergence of $u_{\ve}$ in $C([0,T]: L^2(B_R))$
together with the weak convergence of $p_\ve ={\cal K_{\ve}}u_{\ve}$ and $\nabla p_{\ve}$. So we can pass to the limit in this term. The conclusion is that we have obtained a weak solution of the initial value problem for the equation we have mentioned, posed in $Q_R$ with zero Dirichlet boundary conditions. The regularity of $u$, ${\cal H} u$ and ${\cal K} u$ is as stated before. We also  have the energy formula
\begin{equation}
\int_{B_R} F_\mu (u(t))\,dx+\delta \int_0^t\int_{B_R} \frac{|\nabla u|^2}{d(u)}\,dx dt+\int_0^t\int_{B_R} |\nabla {\cal  H} u|^2\,dx=\int F_\mu(u_0)\,dx.
\end{equation}
\noindent  {\bf Remark.} If we pass also to the limit $\delta\to 0$, which is feasible with no extra effort, then we lose the $H^1$ estimates for $u$, and besides we have a problem with the boundary data that is maybe important. Therefore, we will keep the term $ \delta \Delta u$ for the moment to avoid the problem.

\medskip

\noindent {\bf \ref{sec.ex1}.3}  We will now try to pass to the limit, either as $\mu\to 0$ or as $R\to\infty$, the order depending on convenience. In that sense we recall the detailed form of the energy identity
\begin{eqnarray*}
& \delta \dint_0^t\dint_{B_R} \frac{|\nabla u|^2}{d(u)}\,dx dt+\dint_0^t\dint_{B_R} |\nabla {\cal  H}u|^2\,dxdt +
\int_{B_R} (u_0-u(t))\,dx\\
& \dint_{B_R} (u(t)+\mu)\log(1+\frac{u(t)}{\mu})\,dx =
\dint_{B_R} (u_0+\mu)\log(1+ \frac{u_0}{\mu})\,dx.
\end{eqnarray*}
(Note: we write $u$ but we should write $u_R$ since the solutions changes here with the radius of the ball). Terms 1, 2 and 3 are positive and behave well in both limits. It remains to examine the behavior of Terms 4 and 5. Both are nonnegative which is good for us in the case of Term 4, but  Term 5 diverges as $\mu\to 0$.  Therefore, our choice is letting  $R\to \infty$.  Now   Term 5 is bounded by  hypothesis (though the bound depends on $\mu$). In the limit we easily get a solution of the problem in the whole space, for the equation
 \begin{equation}\label{eq.ws.d}
u_t= \delta \Delta u + \nabla\cdot((u+\mu)\nabla {\cal K}(u)), \qquad x\in \RR^n, \ t>0\,.
\end{equation}
The limit $\delta \to 0$ offers then no difficulty if one wants to take it. However,  the limit $\mu\to0$ needs some extra properties.

\medskip

\noindent {\bf \ref{sec.ex1}.4.} One of such useful properties is the {\sl Conservation of Mass,} that we establish next for the solutions of \eqref{eq.ws.d}.

\begin{lem} \label{le.cm} Under the assumption that $u_0\in L^1(\ren)\cap L^\infty(\ren)$ the constructed nonnegative solutions of the previous problem satisfy
\begin{equation}\label{consmass}
\int u(x,t)\,dx=\int u_0(x)\,dx \qquad \mbox{for all $t>0$.}
\end{equation}
\end{lem}

\noindent {\sc Proof.} Recall that we assume $s<1/2$ if $n=1$. We integrate against a cutoff function to get $\varphi\in C^\infty(\ren)$ supported in $R\le |x|\le 2R$ with $\varphi=1$ for $|x|\le R$. We get
$$
\int u_t\varphi\, dx= \delta \int u\,\Delta \varphi\,dx - \int (u+\mu) (\nabla {\cal K}u\cdot \nabla \varphi)\,dx =I_1+I_2.
$$
For the typical cutoff choice we estimate the first integral as $I_1=O(R^{-2})$ using the fact that $u(t)\in L^1(\ren)$ and then $I_1\to 0$ as $R\to\infty$. As for the last integral, we do
$$
I_2=\int {\cal K}u\, \nabla u\nabla \varphi\,dx+ \int u{\cal K}u\,\Delta \varphi\,dx+
 \mu \int {\cal K}u\,\Delta \varphi\,dx\,.
$$
The latter integral can be estimated as
$$
I_{23}=\mu\int u \,({\cal K}\,\Delta \varphi)\,dx=\mu\|u\|_1O(R^{-2(1-s)})
$$
(where we use the fact that ${\cal K}\Delta$ has homogeneity $2(1-s)$ as a differential operator) and this goes to zero as $R\to\infty$.

Before we estimate the other two integrals we split the kernel $\cal K$ into a bounded part $K_1=\min\{1,{\cal K_s}\}$ and the rest $K_2={\cal K_s}-K_1$ which is supported in a small ball. Both parts are nonnegative and moreover $K_1\in L^\infty$ and $K_2\in L^1$. It means that $K_1*u\in L^q*L^p$ for $q>q_0=n/(n-2s)$  (recall that we always have  $2s<n$) and for every $p\ge 1$, hence
$K_1u\in L^p$ for all $p>q_0$,  while $K_2*u\in L^q*L^p $ with $q<q_0$ and $p\ge 1$ which means that $K_2*u\in L^p$ for all $p\ge 1$. We conclude that ${\cal K}u\in L^p$ for $p>q_0$. We then have
$$
I_{22}= \int u {\cal K}u\,\Delta \varphi\,dx\le \frac{C}{R^2}\|u\|_q\|{\cal K}u\|_p
$$
which has in particular a bound of the form $|I_{22}|\le C\|u\|_1R^{-2}$ that goes to zero as $R\to\infty$. Finally,
$$
I_{21}=\int  {\cal K}u\,(\nabla u\cdot \nabla \varphi)\,dx.
$$
Since $\nabla u\in L^2$ and $\nabla \varphi=O(R^{-1})$ and $\nabla \varphi\in L^{p}$ with $p>n$  we only need
${\cal K}u\in  L^{q}$ with a small $q<2n/(n-2)$ which is true since $q_0<2n/(n-2)$, i.e., $4s<n+2$.

In the limit $R\to\infty$ when $\varphi= 1$, we get (\ref{consmass}).\qed

\begin{thm}\label{thm.ex.app} Let $u_0\in L^1(\ren)\cap L^\infty(\ren)$, $u_0\ge 0$.
Then there exists a weak solution $u=u_{\delta,\ve,\mu}$ of the approximate equation \eqref{appr.eq} posed in $Q_T$ with initial data $u_0$, and \ $u\in L^\infty(0,\infty: L^1(\ren))$, \ $u\in L^\infty(Q)$, \ $\nabla {\cal H}(u)\in L^2(Q)$. Moreover, for all $t>0$ we have
\begin{equation}
\int_{\ren} u(x,t)\,dx=\int_{\ren} u_0(x)\,dx,
\end{equation}
and $\|u(t)\|_\infty \le \|u_0\|_\infty$. The first energy inequality holds,
in the form
\begin{eqnarray*}
& \delta \dint_0^t\dint_{\ren} \frac{|\nabla u|^2}{u+\mu}\,dx dt+\dint_0^t\dint_{\ren} |\nabla {\cal  H}u|^2\,dxdt +\\
& \dint_{\ren} u(t)\log(u(t)+\mu)\,dx + \mu\int_{\ren}\log(1+(u/\mu))\,dx\le\\
&\dint_{\ren} u_0\log(u_0+\mu)\,dx+  \mu\int_{\ren}\log(1+(u_0/\mu))\,dx.
\end{eqnarray*}
\end{thm}
Note that the last integral is less than $\int u_0\,dx$. The parameters $\delta,\ve$ and $\mu$ are larger than 0. Passing to the limits $\delta\to 0$, $\ve\to 0$ offers now no difficulties. But we recall that before taking those limits the solution is smooth in $x$ and $t$, and this may be convenient in justifying calculations to be done below.

\medskip

\section{Tail control. Existence of weak  solutions}
\label{sec.ex2}

We start the section by stating the main existence theorem that will be proved here by passage to the limit $\mu\to 0$ in the construction of the last section. The limit that we want to take  is not trivial because of the possible behaviour of the constructed solutions for large $|x|$ which affects the convergence of some of the resulting integrals. We will pass to the limit in the last formula for the energy  using as an extra tool  a nice decay at infinity, so that the term $\int u\log^{-}(u+\mu)\,dx $ will be bounded uniformly as $\mu \to 0$.

\begin{thm}\label{thm.ex} Let $u_0\in L^1(\ren)\cap L^\infty(\ren)$ be  such that
\begin{equation}
0\le u_0(x)\le A\,e^{-a|x|} \qquad \mbox{for some $A,a>0$}\,.
\end{equation}
Then there exists a weak solution $u$ of Equation \eqref{eq} with initial data $u_0$. Besides,
$u\in C([0,\infty): L^1(\ren))$, $u\in L^\infty(Q)$,  $\nabla {\cal H}(u)\in L^2(Q)$. For all $t>0$ we have
\begin{equation}
\int_{\ren} u(x,t)\,dx=\int_{\ren} u_0(x)\,dx,
\end{equation}
and $\|u(t)\|_\infty \le \|u_0\|_\infty$. The solution decays exponentially as $|x|\to\infty$ as explained in the next tail control subsection. The first energy inequality holds in the form
\begin{equation}
\dint_0^t\dint_{\ren} |\nabla {\cal  H} u|^2\,dxdt +  \dint_{\ren}  u(t)\log(u(t))\,dx
 \le \dint_{\ren} u_0\log(u_0)\,dx\,,
\end{equation}
while the second says that for all $0<t_1<t_2<\infty$
\begin{equation}
\int_{t_1}^{t_2}\int_{\ren} u\,|\nabla {\cal K}u|^2\,dxdt+ \frac 12\int_{\ren} |{\cal H}u(t_2)|^2\,dx\le
\frac 12\int_{\ren} |{\cal H}(u(t_1)|^2\,dx\,.
\end{equation}
\end{thm}

We recall that $Q=\ren\times (0,\infty)$, and we take  $0<s<1$ for all $n\ge 2$, $0<s<1/2$ for $n=1$. The result is first proved for the solutions of the approximate equation with $\delta,\ve,\mu>0$  constructed in Theorem \ref{thm.ex.app}, and the information is then transferred to the original equation once we show that we can pass to the limit and the estimates are uniform.

\subsection { Tail control }

In order to prove the above existence theorem we need to estimate a certain rate of decay of our solutions as $|x|\to\infty$. This will of course depend of a similar assumption that we are imposing on the initial data.  In the case of the PME a possible proof proceeds by constructing explicit weak solutions (or supersolutions) with that property and then using the comparison principle, that holds for that equation. Since we do not have such a general comparison principle here, we have to devise a comparison method with a suitable family of barrier functions that work because they are some kind of  ``exaggerated supersolutions''. What we need is to make sure that {they do not admit a first contact from below.} We will give to functions with such a property the more serious name of {\sl true supersolutions}. This original technique has to be adapted to the peculiar form of the integral kernels involved in operator ${\cal K}$. We prove two kinds of results, the stronger one for small $s$.

\begin{thm}\label{thmex1} Let $0<s<1/2$ and assume that our solution $u$ is bounded $0\le u(x,t)\le L$, and  $u_0$ lies  below a  function of  the form
\begin{equation}
U_0(x)=Ae^{-a|x|}, \ A,a>0.
\end{equation}
If $A$ is large then there is a constant $C>0$ that depends only on $(n,s,a, L, A)$ such that for any $T>0$ we will have the comparison
\begin{equation}\label{compar.barr}
u(x,t)\le Ae^{Ct-a|x|} \quad \mbox{for all $x\in \ren$ and all $0<t\le T$.}
\end{equation}
\end{thm}

 \noindent {\sl Proof. } In order to have enough regularity in the comparison argument below  we make the proof for the solutions constructed in Theorem \ref{thm.ex.app} in the whole space  with parameters $\delta, \ve$ and $\mu>0$ and we will show that the constants in the upper estimate are uniform with respect to such parameters if the
 mentioned  parameters are small.

 \medskip

\noindent $\bullet$  {\bf Reduction.}  By scaling we may put $a=L=1$. This is done by considering instead of $u$ the function $\wu$
defined as
\begin{equation}
u(x,t)= L \,\wu(ax, bt), \quad b=La^{2-2s},
\end{equation}
which satisfies the equation $\wu_t=\delta_1\Delta \wu + \nabla\cdot (\widetilde d(\wu) \nabla K(\wu))$ with $\delta_1=a^{2s}\delta/L.$
Note that then $\wu(x,0)\le A_1\,e^{-|x|}$ with $A_1=A/L$.  A simple calculation then
shows that $C(a,L,A)=La^{2-2s}C(1,1, A/L)$.

\medskip

\noindent $\bullet$ {\bf Contact analysis.}  Therefore, we assume that $0\le u(x,0)\le 1$ and also that
 $$
 u(x,0)\le Ae^{-r} \qquad r=|x|>0,
  $$
 where $A>0$  is a constant that will be chosen below, say larger than 2. Given constants $C,\ve$ and $\eta>0$ we consider  a radially symmetric  candidate to upper barrier function of the form
\begin{equation}
{\widehat U}(x,t)=Ae^{Ct-r}+\ve A\, e^{\eta t},
\end{equation}
 and we take $\ve$ small. Then $C$ will have to be determined in terms of $A$ to satisfy a ``true  supersolution condition'' which is obtained by contradiction at the first point $(x_c,t_c)\in Q_T$ of possible contact of $u$ and ${\widehat U}$.  Note that  if there is contact  it cannot happen only at $|x|=\infty$ since $u$ is integrable and ${\widehat U}$ converges to a positive constant. This is the reason to add the correction term in $\ve$ (use of the correction term can  be avoided if we start the comparison argument with the solutions of the approximate problem posed in a ball $B_R(0)$ with zero boundary data and then pass to the limit $R\to\infty$ as done in the previous section). The contact does not happen at $x_c=0$ if $A$ is large since, putting $r_c=|x_c|$,  we have at the contact point  equality ${\widehat U}=A\,e^{Ct_c-r_c}+\ve Ae^{\eta t_c}=u$ and also $u\le 1$, so that $A\le A\,e^{Ct_c}\le  e^{r_c}$ so $r_c$ is big if  $A$ is. We need at least $A>1$. Note that we are assuming $0<t_c\le T$ and $T$ fixed (in this argument).

 At the point and time of contact we have \ $u=A\,e^{Ct_c-r_c}+\ve A\,e^{\eta t_c}$. Assuming also  that $u$ is $C^2$ smooth, a standard argument also gives
$$
\partial_r u=-A\,e^{Ct_c-r_c}, \quad \Delta u \le A\,e^{Ct_c-r_c}, \quad u_t\ge AC\,e^{Ct_c-r_c}+\ve\eta A\,e^{\eta t_c}\,,
$$
(all of them computed at the point $(x_c,t_c)$), and the spatial derivatives of $u$ at $x_c$ in directions perpendicular to the radius are zero. Next, we put $p={\cal K}u$ and use the equation in full approximate form : \begin{equation}\label{eq.press}
u_t=\delta \Delta u + \nabla u\cdot \nabla p+(u +\mu)\Delta p,
\end{equation}
 to get the basic inequality
$$
CA\,e^{Ct_c-r_c}+\ve A\eta\,e^{\eta t_c} \le \delta A\,e^{Ct_c-r_c} -A\,e^{Ct_c-r_c}\overline {\partial_r p}+ (A\,e^{Ct_c-r_c}+\ve A\,e^{\eta t_c}+\mu) \overline {\Delta p}.
$$
Cleaning up this expression, we get at the contact point the following condition:
\begin{equation}\label{ineq.n}
C+ \ve \eta e^{\xi} \le \delta  -\,\overline {\partial_r p}+ (1+ \ve\, e^{\xi}+\mu_1) \overline {\Delta p}\,,
\end{equation}
where $\xi=r_c+(\eta-C)t_c$, $\mu_1=(\mu/A)\,e^{r_c-Ct_c}$ and the overline indicates that  the values of $\partial_r p$ and $\Delta p$
are calculated at the point of contact.

\medskip

\noindent {\bf Summing up the main ideas.} In \eqref{eq.press} we have written the nonlinear term of the equation involving the fractional operator in ``non-divergence form'', precisely as the sum of a transport term involving first derivatives on $p$ and another term containing the Laplacian of $p$. When we evaluate those terms at the contact point with the barrier and arrive at \eqref{ineq.n}, the latter term is not difficult in terms of global integrals of $u$ and the other two first terms in \eqref{eq.press} have also simple contributions. However, the transport term does not have a sign and its control becomes more difficult. We proceed with the simple cases and
leave the difficult case for the end.

\medskip

$\bullet$  In order to get a contradiction with  inequality \eqref{ineq.n} we will
 estimate the values of $\overline {\partial_r p}$ and $\overline {\Delta  p}.$
  For $n\ge 1$, $0<s< 1$, we use the formula
\begin{equation}\label{eq.tc2.n}
p(x,t)=c\int \frac{u(x-y,t)}{|y|^{n-2s}}dy=c\int
\frac{u(x+y,t)}{|y|^{n-2s}}dy,
\end{equation}

with some $c=c(s)>0$, which produces the singular integrals
\begin{equation}\label{eq.tc3.n2}
p_{x_i}(x,t)=c_1\int \frac{(u(x+y,t)-u(x-y,t))\,y_i}{2|y|^{n+2-2s}}\,dy,
\end{equation}
($i=1,\cdots, n$) that for $s<1/2$ can also be written as
\begin{equation}\label{eq.tc3.n}
p_{x_i}(x,t)=c_1\int \frac{(u(x+y,t)-u(x,t))\,y_i}{|y|^{n+2-2s}}\,dy,
\end{equation}
and finally,
\begin{equation}\label{eq.tc4.n}
\Delta p(x,t)=c_2\int
\frac{u(x+y,t)+u(x-y,t)-2u(x,t)}{|y|^{n+2-2s}}\,dy.
\end{equation}
 Here and in the sequel of the proof we denote  by $c,c_i, K, K_i$ different absolute positive constants, i.,e., constants that depend only on $n$ and $s$.

\medskip

$\bullet$ We now estimate $\overline {\Delta p}$ at $x_c$. In view
of inequality (\ref{ineq.n}) we need a control from above. Using
formula \eqref{eq.tc4.n} we see that the integral on the set $|y|\ge 1$ is
absolutely convergent and  absolutely bounded  (since $u$ is bounded by 1).
Besides, to bound from above the part of integral on the set  $|y|\le 1$  we will
use the fact that $u(t_c)$ lies below $Ae^{Ct_c-|x|}$ with tangency at $|x|=r_c$; since
the numerator of the integrand is the discretization of the second derivative for $u$
at $x_c$ we use the fact that it is bounded above by the same expression for $\widehat U$ to get
$$
\int_{|y|\le 1}
\frac{u(x_c+y,t_c)+u(x-y,t_c)-2u(x,t_c)}{|y|^{n+2-2s}}\,dy\le \int_{|y|\le
1} \frac{2Ae^{Ct_c-r_c}|y|^2}{|y|^{n+2-2s}}\,dy,
$$
which is bounded since $s>0$. Noting that $Ae^{Ct_c-r_c}\le 1$ we conclude that the term in
(\ref{ineq.n}) containing $\overline {\Delta p}$  contributes with
at most $K_1$ to the inequality.

$\bullet$ Next, we estimate the transport term involving $\overline {\partial_r p}$. We have
$$
\overline {\partial_r p} =c_1 \int \frac{(u(x+y,t)-u(x,t))(\hat x_c\cdot \hat y)}{|y|^{n+1-2s}}\,dy,
$$
where $\hat x=x/|x|$ and $\hat y=y/|y|$. Since we want to
get a contradiction in formula \eqref{ineq.n} and  $\overline {\partial_r p}$ appears with
a minus sign, we need to estimate this term from below. This integral is delicate
so we split the calculation  into several pieces.  At this moment we make the further assumption
 $s<1/2$. Then the integral for $|y|\ge 1/2$ is bounded as:
$$
|I_{\{|y|>1/2\}}(\overline {\partial_r p})|\le  c\displaystyle\int_{|y|>1/2}
 \frac{|u(x_c+y,t_c)-u(x_c-y,t_c)|}{|y|^{n+1-2s}}\,dy\le c\int_{y\ge 1}
\frac{1}{|y|^{n+1-2s}}\,dy\le K_2.
$$
Hence, this part has the desired control.

\noindent $\bullet$  The integral on the ball $\{|y|\le 1/2\}$ is split again into two parts. By rotation we may assume that $x_c$ is directed along the first axis, $x_c=(r, 0,\cdots, 0)$, where $r=|x|$. Then the integral is calculated on $\Omega_1=\{y: |y|\le  1/2, y_1>0\}$ and on $\Omega_2= \{y: |y|\le  1/2, y_1<0\}$. The last is easily bounded since $u$ touches $\overline U$ at $x_c$ and lies below everywhere at  time $t_c$. Hence,
$$
\begin{array}{c}
-I_{\Omega_2}(\overline {\partial_r p})=\displaystyle\int_{\Omega_2}
 \frac{(u(x_c,t_c)-u(x_c+y,t_c))y_1}{|y|^{n+2-2s}}\,dy=
\displaystyle\int_{\Omega_1}\frac{(u(x_c-y,t_c)-u(x_c, t_c))y_1}{|y|^{n+2-2s}}\,dy\\
    \le \displaystyle\int_{\Omega_1}
\frac{(\overline{U}(x_c-y,t_c)-\overline{U}(x_c, t_c))y_1}{|y|^{n+2-2s}}\,dy\le
\displaystyle\int_{\Omega_1}
\frac{2Ae^{Ct_c-r_c}}{|y|^{n+1-2s}}\,dy\le K_3.
\end{array}
$$
where in the last inequality we have used the linear approximation for $\overline{U}$. Again, this is a good control.

\noindent $\bullet$  Now, the difficult part. The integral on $\Omega_1$ (i.\,e., the ``half integral looking outside near $x_c$'') is more delicate since the difference $u(x,t_c)-u(x+y,t_c)$ could in principle drop quite abruptly even at a relatively short distance from $x_c$ and this would make the integral very big. However, we have a miraculous control by combining the estimate on the Laplacian with the good part of the estimate of $\overline {\partial_r p}$.

\begin{lem} With the previous assumptions and notations we have
\begin{equation}\label{miracle}
-I_{\Omega_1}(\overline {\partial_r p})+\frac12 I_{\Omega_1}(\overline {\Delta p})\le -\frac12 I_{\Omega_2}(\overline {\partial_r p}).
\end{equation}
\end{lem}

\noindent {\sl Proof.} We combine the integral of $-\overline {\partial_r p} $ with a part of the integral for $\overline {\Delta p}$ as follows:
$$
Y= -\displaystyle \int_{\Omega_1}  \frac{(u(x_c+y,t_c)-u(x_c,t_c))y_1}{|y|^{n+2-2s}}+ \frac12 \int_{\Omega_1}
\frac{u(x_c+y,t_c)+u(x_c-y,t_c)-2u(x_c,t_c)}{|y|^{n+2-2s}}\,dy
$$
and we study carefully the integrand of $Y$. We have a numerator of the form
$$
\begin{array}{c}
\displaystyle(u(x_c,t_c)-u(x_c+y,t_c))y_1 + \frac12(u(x_c+y,t_c)+u(x_c-y,t_c)-2u(x_c,t_c))=\\
\displaystyle-(\frac12-y_1)(u(x_c,t_c)-u(x_c+y,t_c))+ \frac12(u(x_c-y,t_c)-u(x_c,t_c))
\end{array}
$$
Luckily,  the first term is negative, hence we conclude that estimate \eqref{miracle} holds. \qed

\medskip

\noindent $\bullet$ We can now finish the contradiction argument. All these estimates allow to conclude that $-I(\overline {\partial_r p})\le K_4$. We go now back to \eqref{ineq.n} and conclude that the inequality implies that
$$
C+ \ve (\eta-K) e^{r_c+(\eta-C)t_c}  \le \delta + K+\frac{K\mu}{A}\,e^{r_c},
$$
This cannot happen if we choose  $\ve\ge \frac{\mu}{A}$ and $C=\eta \ge 2K$, where we keep the extra condition that $\mu$ and $\delta$ must be small (at least less than 1).  Then  $\ve$ may go to zero as $\mu\to 0$.

\medskip

$\bullet$ Finally, we worry about the regularity of the
solutions. To make the above proof fully rigorous, we apply the argument to the solutions of the regularized problem where we smooth the velocity field $\nabla p$ by regularizing the kernel. We have to make the estimates for ${\overline {p}_{r}}$ and $\overline {\Delta p}$ when
$$
p(x,t)=\int K_\ve(y)u(x+y)\,dy.
$$
where for instance $K_\ve(y)=K(y)$ for $\ve\le |y|\le 1/ \ve$,
$K_\ve(y)$ is a parabolic cap with $C^1$ fit in $|y|\le
\ve$, and finally $K_\ve(y)=0$ for $|y|\le 2/\ve$. The regularization mentioned in Section \ref{sec.ex1} will also do. The
solutions $u_\ve$ to this problem have bounded speeds, they are
smooth and bounded with smooth and bounded $p_x$ and the previous
estimates for ${\overline p}_{r}$ and $\overline{\Delta p}$  hold
uniformly in $\ve$. Passing to the limits $\ve\to 0$, the previous conclusions hold for any weak limit solution
as constructed above, cf. equation \eqref{eq.ws.d}. The extra limit $\delta \to 0$ offers then no difficulty. \qed

\begin{thm}\label{thmex2} Let now $1/2\le s<1$. Under the assumptions of the previous theorem the stated tail estimate works locally in time. The global statement must be replaced by the following: there exists an increasing function
$C(t)$ such that
\begin{equation}\label{compar.barr.sg12}
u(x,t)\le Ae^{C(t)t-a|x|} \quad \mbox{for all $x\in \ren$ and all $0<t\le T$.}
\end{equation}
\end{thm}

\noindent {\sl Proof.} (1) In the previous proof we had to put $s<1/2$ only because of the problem in estimating the integral for $\partial_r p$ on an exterior domain, away from the contact point. When $s\ge 1/2$ we can estimate such integral as a convolution integral between $u$ and the kernel $K_1(y)=y_1|y|^{-n-2+2s}\chi(|y|\ge 1)$. Now,  this kernel belongs to $L^p(\ren)$ for all $p>n/(n+1-2s)$, hence we only need to bound $u(t)$ in an $L^q(\ren)$ norm with $1\le q<n/(2s-1)$ to get
\begin{equation}
I_*:=|I_{\{|y|>1/2\}}(\overline {\partial_r p})|\le \|u(t)\|_q\|K_1\|_p
\end{equation}
Moreover, since $u\in L^1(\ren)\cap L^\infty(\ren)$ we have
$$
\|u(t)\|_q\le \|u(t)\|_1^{1/q}\|u(t)\|_\infty^{(q-1)/q}
$$
We know that $\|u(t)\|_\infty\le 1$ and $u(x,t)\le Ae^{Ct}e^{-|x|}$ hence
$\|u(t)\|_1\le cAe^{Ct}$. Therefore the  term contributes
$$
I_* \le KA^{1/q}e^{Ct_c/q}
$$
which allows to go back to \eqref{ineq.n} and get
\begin{equation}\label{ineq.sg12}
C+ \ve (\eta-K) e^{r_c+(\eta-C)t_c}  \le \delta + K+ KA^{1/q}e^{Ct_c/q}+ \frac{K\mu}{A}\,e^{r_c},
\end{equation}
The contradiction argument works as before with only one big difference. Once we put $C=\eta=KA^{1/q}$ the contradiction is gotten if we restrict the time so that $e^{Ct_c/q}\le 2$
which happens if
$$
t_c\le T_1= c_2/C=c_3K^{-1}A^{-1/q}= c_4A^{-1/q}. $$
We do not play with $A$ here, only $A$ bigger than 2 or the like. 

\medskip

(2) Once we prove estimate (\ref{compar.barr}) for $0<t<T_1$ we can repeat the argument for
another time interval, but now with constant $A_1=Ae^{CT_1}= Ae^{c_2}$, and then $C_1=KA_1^{1/q}$. We get a valid time interval
$$
T_2=c_4A_1^{-1/q}=c_4 A^{-1/q}e^{-c_2/q}
$$
where a new factor will appear every time we repeat the iterations. In this way we can extend the bound up
to a certain time that depends on the initial data through the value of $A$.

\medskip

(3) In order to prove that the time for which the estimate is valid goes forever we need to improve the $L^1$ norm of $u(t)$ by using the fact that $u\le 1$ together with $u(x,t)\le Ae^{Ct}e^{-|x|}$. Summing the contributions of both upper bounds, we now get
$$
\|u(t)\|_1 \le c_0 (\log (Ae^{Ct}))^n\left(1+ A^{-1}e^{-Ct}\right)
$$
so that
$$
I_* \le K (\log A+ Ct)^{n/q}\left(1+ A^{-1}e^{-Ct}\right)^{1/q}
$$
and then Inequality \eqref{ineq.sg12} may be replaced by
\begin{equation}\label{ineq.sg12.2}
C+ ...  \le  K+ K(\log A+Ct)^{n/q}\left(1+ A^{-1}e^{-Ct}\right)^{1/q} + \frac{K\mu}{A}\,e^{r_c},
\end{equation}
where we have dropped the terms in $\delta, \ve$ and $\mu$ that add no novelties or problems. Now we may put $C=K[(\log A)^{n/q}+2]$ and $CT_1=c_2$ so that
$$
T_1=c_4(\log A)^{-n/q}
$$
At the new starting time, we have $A_1=Ae^{CT_1}=Ae^{c_2}$ and then
$$
T_2=c_4(\log A_1)^{-n/q}=c_4(\log A +c_2)^{-n/q}
$$
and so on. Since for large $k$ we get $T_k\sim c k^{-n/q} $ and we m ay always take $q>n$, the series $\sum T_k$ diverges, so that Estimate \eqref{compar.barr.sg12} is global in time.  \qed

\subsection{Proof of the existence result}

We may now pass to the limit in the weak formulation of the equation satisfied by $u=u_{\delta,\ve, \mu}$
constructed at the end of previous section. The exponential decay bound on the solutions, which
is uniform in $\mu$ allows to improve the consequences of the energy inequality that is now
written as
\begin{eqnarray*}
& \delta \dint_0^t\dint \frac{|\nabla u_\mu|^2}{u_\mu +\mu}\,dx dt+\dint_0^t\dint |\nabla {\cal  H}u_\mu|^2\,dxdt +  \dint u_\mu(t)\log^+(u_\mu(t)+\mu)\,dx  \\
& \le
\dint u_0\log(u_0+\mu)\,dx+\mu\dint \log(1+(u_0/\mu))\,dx+\dint u_\mu(t)\log^-(u_\mu(t)+\mu)\,dx .
\end{eqnarray*}
where we use the notation $u_\mu$ for the solution for clarity. Since the right-hand side is bounded uniformly in $\mu$ we get a uniform estimate for the family
$|\nabla {\cal  H}u_\mu|$ in $L^2(Q)$. We also have a uniform bound on $\nabla u_\mu$ in $L^2(Q)$ if $\delta>0$. We can therefore pass to the weak limit in $u_\mu \to \tilde u$, and then ${\cal  H}u_\mu\to {\cal  H}\tilde  u$ and $\nabla {\cal  H}u_\mu \to \nabla {\cal  H}\tilde  u$. The same happens to ${\cal  K}u_\mu$ and $\nabla {\cal  K}u_\mu$. The weak solution is obtained and we have
 \begin{eqnarray*}
& \delta \dint_0^t\dint \frac{|\nabla \tilde  u|^2}{\tilde u}\,dx dt+\dint_0^t\dint |\nabla {\cal  H}\tilde u|^2\,dxdt +  \dint \tilde u(t)\log^+(\tilde u(t))\,dx  \\
& \le
\dint u_0\log(u_0)\,dx+\dint \tilde u(t)\log^-(\tilde u(t))\,dx .
\end{eqnarray*}
Note that the term $\int \mu \log (1+(u_0/\mu))\,dx$ disappears in the limit by the Dominated Convergence Theorem since the integrand is uniformly bounded by $u_0$. Due to the decay at infinity the integral $\int \tilde u(t)\log \tilde u(t)\,dx$
is absolutely convergent.

The constructed solution has the same $L^\infty$ bound as $u_0$ and conservation of mass holds by
the proof of Lemma \ref{le.cm}.

Passing to the limit $\delta\to 0$ and $\ve\to 0$ offers no difficulty, so Theorem \ref{thm.ex}
is proved but for the second energy estimate. \qed

\subsection{Second energy estimate}

Let us establish the second energy estimate for weak solutions with tail decay. We compute formally
of the approximations where everything is justified
\begin{equation*}
\begin{array}{l}
\dfrac 12\frac{d}{dt}\dint \varphi|{\cal H}u(x,t)|^2\,dx= \int \varphi \,{\cal H}u { \cal H}u_t\,dx= \dint \,{\cal H}(\varphi {\cal H}u)\,u_t\,dx= \\ \dint {\cal H}(\varphi {\cal H}u)\, \nabla (u\nabla {\cal K}u)= -\dint u\,\nabla {\cal H}(\varphi {\cal H}u)\cdot \nabla {\cal K}u\,dx=\\
-\dint u\,{\cal H}(\varphi \nabla {\cal H}u)\cdot \nabla {\cal K}u\,dx-
\dint u\,{\cal H}(\nabla \varphi  {\cal H}u)\cdot \nabla {\cal K}u\,dx,
\end{array}
\end{equation*}
so that, putting $\psi=1-\varphi$ we have
\begin{equation*}
\begin{array}{l}
\dfrac 12\frac{d}{dt}\dint \varphi|{\cal H}u(x,t)|^2\,dx +\int \varphi u |\nabla {\cal K}u|^2\,dx=\\
\dint u\,{(\psi \nabla {\cal K}u-\cal H}(\psi \nabla {\cal H}u))\cdot \nabla {\cal K}u\,dx-
2\dint   {\cal H}u \,\nabla \varphi \cdot {\cal H}(u\, \nabla {\cal K}u)\,dx,
\end {array}
\end{equation*}
With the properties we have for $u$ the right-hand side must go to zero as $\varphi\to 1$
along the typical cutoff sequence, at least after integration in time.

\section{Finite propagation. Solutions with compact support}
\label{sec.fp}

One of the most important features of the porous medium equation
and other related degenerate parabolic equations is the property
of finite propagation, whereby  compactly supported initial data
$u_0(x)$ give rise to solutions $u(x,t)$ that have the same
property for all positive times, i.e., the support of $u(\cdot,t)$
is contained in a ball $B_{R(t)}(0)$ for all $t>0$ and $R(t)$ is bounded on bounded intervals $0<t<T$.

A possible proof in the case of the PME proceeds by constructing explicit weak solutions with that property (and possibly larger  initial data) and then using the comparison principle, which holds for that equation. Since we do not have such a general principle here, we have to devise a comparison method with a suitable family of upper barriers that behave as ``true (exaggerate) supersolutions''. The technique has been presented in whole detail in the tail analysis of the previous section and is here adapted to the peculiar needs of bounded support. Here is the end result.

\begin{thm}\label{prop.sc1} Assume that $u$ is a bounded solution, $0\le u\le L$, of equation \eqref{eq} with \  ${\cal K}=(-\Delta)^{-s}$ with $0<s<1$ $(0<s<1/2$ if $n=1)$, as constructed in Theorems {\rm \ref{thmex1}} and {\rm \ref{thmex2}}.  Assume that $u_0$   has compact support. Then $u(\cdot,t)$ is compactly supported for all $t>0$. More precisely, if  $0<s<1/2$ and $u_0$ is below the ''parabola-like'' function
\begin{equation}
U_0(x)=a(|x|-b)^2,
\end{equation}
for some $a,b>0$, with support in the ball $B_b(0)$, then there is $C$  large enough, such that
\begin{equation}
u(x,t)\le a(Ct-(|x|-b))^2
\end{equation}
 Actually, we can take $C=C_0(n,s)L^{(1/2)+s}a^{(1/2)-s}$. For $1/2\le s<2$ a similar conclusion is true, but now $C$ is an increasing function of $t$ and we do not obtain a scaling formula for its dependence of $L$ and $a$.
\end{thm}

\noindent {\sl Proof.}
The application of the  method is very similar to the case worked out  the tail control section. Therefore, we will dispense with some of the technicalities of regularization to gain space and clarity.  We assume that our solution $u(x,t)\ge 0$ has bounded initial data $u_0(x)=u(x,t_0)\le L$ and also  that $u_0$ is below the parabola \
$ U_0(x)=a(|x|-b)^2, \ a,b>0; $ moreover,  the support of $u_0$ in the ball of radius $b$ and
the graphs of $u_0$ and $U_0$ are strictly separated in that ball. We take as comparison function \ $U(x,t)=a(Ct-(|x|-b))^2$ and argue at the first point and time where $u(x,t)$ touches $U$ from below. The fact that such a first contact point happens for $t>0$ and does not happen at $x=\infty$ is justified by regularization as before. We put $r=|x|$.

By scaling we may put $a=L=1$. See detail of the reduction step below.
 We examine in detail the situation in which the touching  point $(x_c,t_c)$ is not the minimum, say,
 $x_c$ lies at a distance from the front $|x_f(t)|:=b+Ct$, so that
  $b+Ct_c-|x_c|=h>0$. Note that since $u\le 1$ we must have
$|h|\le 1$. Assuming also  that $u$ is $C^2$ smooth, a standard
argument gives
$$
u=h^2, \quad u_r=-2h, \quad \Delta u \le 2n, \quad u_t\ge 2Ch,
$$
all of them computed at the point $(x,t)=(x_c,t_c)$. Putting $p={\cal K}u$ and using the
equation $u_t=\nabla u\cdot \nabla p +u\Delta p$, we get the inequality
\begin{equation}\label{eq.cs1}
2Ch\le -2h{\overline {p_r}}+ h^2{\overline {\Delta p}},
\end{equation}
where the overline indicates that the values of $p_r$ and $\Delta p$
are calculated at the point of contact. Moreover, $u(x,t_c)\le (x_f-x)^2$
for all $x\in \ren$. In order to get a contradiction we will
 estimate the values of ${\overline p}_{r}$ and $\overline {\Delta p}$. The formulas for $p$, $p_r$ and $\Delta p$ in terms of $u$ are given
in \eqref{eq.tc2.n}, \eqref{eq.tc3.n2}, \eqref{eq.tc3.n}, \eqref{eq.tc4.n}.

$\bullet$ Estimating ${\overline {\Delta p}}$ at the contact point offers no novelties. As  before and in view of inequality (\ref{eq.cs1}) we need a control from above. Using
formula \eqref{eq.tc4.n} we see that the integral on $|y|\ge 1$ is
absolutely convergent and  bounded (since $u$ is bounded).
Besides, to bound the integral for $|y|\le 1$ from above we will
use the fact that $u$ lies below a parabola, hence
$$
\int_{|y|\le 1}
\frac{u(x+y,t_c)+u(x-y,t_c)-2u(x,t_c)}{|y|^{n+2-2s}}\,dy\le \int_{|y|\le
1} \frac{2|y|^2}{|y|^{n+2-2s}}\,dy,
$$
which is bounded since $s>0$. We conclude that the term in
(\ref{eq.cs1}) containing ${\overline p}_{xx}$  contributes with
at most $Kh^2$ to the inequality, where $K>0$ is an absolute
constant.

$\bullet$ Next, we estimate ${\overline {p}_r}$. Much of the argument is similar to the tail
analysis but the end is more delicate. Since we want to
get a contradiction in formula \eqref{eq.cs1} and  the coefficient
$u_r =-2h$ of ${\overline p}_{r}$ is negative, we need to estimate
this term from below. Using formula \eqref{eq.tc3.n2} we see again that the integral for $|y|\ge 1/2$
is absolutely and uniformly bounded by a constant $K_2$.

\noindent $\bullet$  The integral on the ball $\{|y|\le 1/2\}$ is split into several  parts.
We will drop for convenience of writing the dependence on $t_c$ in the formulas the follow.  By rotation we may assume that $x_c$ is directed along the first axis, $x_c=(r, 0,\cdots, 0)$, where $r=|x|$. Then the integral is calculated on $\Omega_1=\{y: |y|\le  1/2, y_1>0\}$ and on $\Omega_2= \{y: |y|\le  1/2, y_1<0\}$. The last is easily bounded since $u$ touches $\overline U$ at $x_c$ and lies below everywhere at  time $t_c$. As in the tail analysis we have,
$$
\begin{array}{c}
-I_{\Omega_2}(\overline {\partial_r p})=\displaystyle\int_{\Omega_2}
 \frac{(u(x_c)-u(x_c+y))y_1}{|y|^{n+2-2s}}\,dy\le K_3\,.
\end{array}
$$

\noindent $\bullet$  The integral on $\Omega_1$ (i.\,e., the ``half integral looking outside near $x_c$'') is more delicate as we have said, since the difference $u(x,t_c)-u(x+y,t_c)$ could in principle drop quite abruptly even at a relatively short distance from $x_c$ and this would make the integral very big.
 There is a part with $|y|$ between $\theta h$ and $1$ for any given $\theta\in (0,1)$ that is easy (note that
$0\le u(x,0)-u(x+y,t)\le u(x,0)=h^2$):
$$
\displaystyle \left| \int_{\theta h}^{1}
\frac{u(x_c+y)-u(x_c)}{|y|^{n+1-2s}}\,dy \right| \le h^2\int_{\theta
h}^{1} \frac{1}{|y|^{n+1-2s}}\,dy =  \ ch^2  (\theta
h)^{-1+2s}=c_3(\theta)h^{1+2s}
$$
 and this is good even for small $h$.

\medskip

$\bullet$  The last part of the integral for ${\overline p}_r$, over the half-ball $H_1=\{0<|y|<\theta h\}$ with $y_1>0$, is in principle bad since \ $u(x_c+y)$ could drop abruptly, thus making the integral very negative or even divergent near $y=0$.  We are going to combine the integral of $-\overline {\partial_r p} $ with a part of the integral for $\overline {\Delta p}$ as follows:
$$
Y= -\displaystyle \int_{H_1}  \frac{(u(x_c+y)-u(x_c)y_1}{|y|^{n+2-2s}}+ \frac{h}2 \int_{H_1}
\frac{u(x_c+y)+u(x_c-y)-2u(x_c)}{|y|^{n+2-2s}}\,dy
$$
and we study carefully the integrand of $Y$. We have a numerator of the form
$$
\begin{array}{c}
\displaystyle(u(x_c)-u(x_c+y)y_1 + \frac{h}{2}(u(x_c+y)+u(x_c-y)-2u(x_c))=\\
\displaystyle-(\frac{h}{2}-y_1)(u(x_c)-u(x_c+y))+ \frac{h}{2}(u(x_c-y)-u(x_c))
\end{array}
$$
Luckily,  the first term is negative (we take $0<\theta <1/2$ (a security factor), hence we conclude that
$$
-I_{H_1}(\overline {\partial_r p})+\frac{h}2  I_{H_1}(\overline {\Delta p})\le -\frac{h}2 I_{\Omega_2}(\overline {\partial_r p}).
$$

$\bullet$  We may now sum up all the terms in the right-hand side of \eqref{eq.cs1} and show that they are bounded above by $Kh$ where $K$ is a uniform constant. Therefore, for large $C$ inequality (\ref{eq.cs1}) is impossible, hence there
cannot be a contact point with $h\ne 0$. In this way we get a
minimal constant $C=C_0(n,s)$ for which such contact does not take
place.

\noindent $\bullet$ {\sc Reduction Step.} We use it to get the dependence on $L$
and $a$. Here, it goes as follows: since the equation scaling is
$$
\widehat u(x,t)=Au(Bx,Tt)
$$
with parameters $A,B,T>0$ such that $T=AB^{2-2s}$, if we have done
the proof for $u$ that has height 1 and is below $(|x|-b_0)^2$
initially, and get a comparison with a $U$ as above with speed
$C_0>0$, then the assumptions $\widehat u(x,t)\le L$ and $\widehat
u(x,t_0)\le a(|x|-b)^2$ initially, are satisfied if we put
$$
A=L, \quad AB^2=a, \quad b=b_0/B,
$$
i.e., $A=L$, $B=(a/L)^{1/2}$, and then $T=a^{1-s}L^s$. The new
speed is then
$$
\widehat C=\frac{C_0T}B=C_0a^{1/2-s}L^{1/2+s}.
$$

\noindent $\bullet$ We still have to consider the modification of the proof when $1/2\le s<1$. The only problem is
the estimate of $\overline{\partial_r p}$ on the exterior of a ball. This is done as in the tail control case.

\medskip

 We next lemmas complete the details of the comparison proof that
have been left out in the previous lemma.

\begin{lem} Under the assumptions of Prop. {\rm \ref{prop.sc1}}
there is no contact either at $u=0$, in the sense that strict separation of $u$ and
$U$ holds for all $t>0$ if $C$ is large enough as in the previous lemma.
\end{lem}

\noindent {\sl Proof.} Precisely, what we want is to eliminate the
possible contact of the supports at the lower part of the
parabola. Instead of doing this  by analyzing the possible contact
point, we proceed by a change in the test function that we replace
by
$$
U_\ve= (Ct-(|x|-b))^2 +\ve(1 + D t) \quad \mbox{for} \quad |x|\le
b+Ct,
$$
and $U_\ve=\ve(1 + D t)$ for $1x|\ge b+Ct$. Here $\ve >0$ is a small
constant and $D>0$ will be suitably chosen. Assume that the
solution starts at $t=0$ and touches  this test function
$U_\ve$ for the first time at the time  $t=t_c$.

The argument for a contact at the points
$|x|< b+Ct_c$ where $U_\ve$ is a parabola works like previously.

The argument at contact points $|x|\ge b+Ct_c$ leads to
\begin{equation}\label{eq.cs1b}
D\ve \le \ve (1+D t_c){\overline {\Delta p}},
\end{equation}
where the overline indicates as before that the value $\Delta p$ are
calculated at the point of contact. Moreover, $u(x,t_c)\le U_1(x,t_c)$
for all $x\in \RR$. If we are able to prove again that ${\overline
 {\Delta p}}\le K$ we will get
$$
D\ve \le \ve (1+Dt_c)K
$$
which is contradictory if $D> K$, say $D=2K$, and $t_c<
1/D=1/(2K)$, This estimate is uniform in $\ve$ and gives in the
limit and upper bound of the form $u(x,t)
\le (Ct-(|x|-b))^2$ for all
$0<t<1/2K$, and as a corollary, the support of $u$  bounded on the right by the
line $|x|=Ct+b$ in that time interval. After this time the process
can be repeated and the conclusion is true for all times.
\qed

\noindent $\bullet$  We now reflect for a moment on the regularity requirements.
Arguing as in the tail control case by using the smooth solutions of the approximate equations, the previous conclusions hold for any weak limit solution.

\subsection{Consequences. Growth estimates of the support}

 The following  analysis is done for $s<1/2$ and concerns  bounded
solutions with compactly supported initial data. By free boundary we mean, as usual, the topological boundary
of the support of the solution ${\cal S}(u):=\overline{\{(x,t): u>0\}}$.

\begin{cor} Let $u_0$ be bounded above by $L$ with $u_0(x)=0$ for $|x|\le R$.  Then, we get an
estimate for the free boundary points of the form $|x(t)|\le R+
C_2\,t^{1/(2-2s)}$ if $s< 1/2$.\end{cor}

\noindent {\sl Proof:} We know that the support of $u(\cdot,t)$ is bounded for all times, say, it is contained in a ball
of radius $r(t)$.  We take a time $t_1$ and find a parabolic barrier as before, with coefficient $a>0$ that is initially above and separated from $u(\cdot,t_1)$. This can be done by choosing first $r_1$ such that  $ar_1^2=L$ and then
putting in the formula of the parabolic barrier  $b=r(t_1)+r_1+ \epsilon$, and then we can go
forever in time in the comparison.   Using the speed estimate in Theorem \ref{prop.sc1} we get
in the limit $\epsilon\to 0$:
$$
r(t)- r(t_1)-r_1\le C(t-t_1)=C_0
L(t_1)^{1/2+s}a^{1/2-s}(t-t_1)=C_0\,L(t_1)(t-t_1)/r_1^{1-2s}º
$$
Here $r_1>0$ is free by moving $a$.  We can use the $L^\infty$ bound $L(t_1)\le L(0)$
We get
$$
r(t)\le r(t_1)+r_1 + C_0\,L\,(t-t_1)/r_1^{1-2s}º
$$
Now take  $t_1=0$ and optimize the right-hand expression  in $r_1>0$
for $s<1/2$. . Notice that  in the limit $s=1/2$ we would get linear growth,
while for $s=0$ we get the standard $t^{1/2}$ growth of the Porous Medium Equations under these assumptions.

\section{Persistence of positivity}
\label{sec.pers} \normalcolor

We establish another property that plays an interesting role in the theory of porous medium equations to avoid degeneracy points for the solutions. It is called persistence of positivity. For continuous solutions it implies non-shrinking of the support.

 \begin{thm} Let $u$ be a weak solution as constructed in Theorem \eqref{thm.ex} and assume that $u_0(x)$ is positive
 in a neighborhood of a point $x_0$. Then $u(x_0,t)$ is positive for all times $t>0$.
 \end{thm}

 \noindent {\sl Proof.} This issue allows for another use of the technique presented in the tail analysis, but this time with a true subsolutions. We assume that $u_0(x)\ge c>0$ in a ball $B_R(x_0)$. By translation and scaling we may assume that $x_0=0$ and $c=R=1$. The  idea is to study the contact
point with a parabola that shrinks quickly in time, like
$$
U(x,t)=e^{-at}F(|x|),
$$
with $F$ suitably chosen and $a>0$ large enough. Firstly, we choose $F$ to be radially symmetric and decreasing,
with $F(0)=1/2$ and $F(|x|)=0$ for $|x|>1/2$.  The contact point $(x_c,t_c) $  is sought in
$B_{1/2}(0)\times (0,\infty)$. By approximation we may assume that $u$ is positive everywhere
so no contact at the parabolic border is assumed. At a positive
contact point we have $u_t\le U_t= -aU$,
$\nabla u=e^{-at}F'(|x|){\bf e_r}$ and the standard arguments on the equation imply
$$
-aF(|x|) \ge F'(|x|)\overline {p_r}+
F(|x|)\overline{\Delta_x p}
$$
Now we have
$$
p(x,t)={\cal K}u(x,t)=e^{-at}({\cal K}F)(r)
$$
As in the tall control analysis we prove that $\overline{\Delta_x p}$ is bounded uniformly.
The novelty is that $F'\le 0$ implies that $p_r\le 0$ so that the term $F'(|x|)\overline {p_r}\ge 0$.
In this way the inequality implies
$$
a\le Ke^{-at}
$$
which is false if $a>K$. Hence, under this size assumption there can be no contact point, and $U$ is a true
subsolution, and positive for all times in $B_{1/2}(0)$. \qed


\section{Appendix. Fractional Laplacians}\label{appendix}

\normalcolor

 We collect some data on fractional Laplacians for the
reader's convenience

\noindent {\bf Fractional Laplacians and potentials.} According to
Stein \cite{Stein}, Chapter V, the definition of
$(-\Delta)^{\beta/2}$ is done by means of Fourier series
$$
  ((-\Delta)^{\beta/2}f)^{\widehat{}}(x)=(2\pi|x|^{\beta} \hat f (x)
$$
and can be used for positive and negative values of $\beta$. For
$\beta=-\alpha$ negative, with $0<\alpha<n$, we have the equivalence
with the Riesz potentials \color{blue} \cite{Riesz} \normalcolor
\begin{equation}
      (-\Delta)^{-\alpha/2}f =I_{\alpha}(f)
      :=\frac1{\gamma(\alpha)}\int_{\RR^n} \frac{f(y)}{|x-y|^{n-\alpha}}dy
\end{equation}
(acting on functions of the class $\cal S$ for instance) with
precise constant
$$
\gamma(\alpha) =\pi^{n/2}2^{\alpha}\Gamma(\alpha/2) /
\Gamma((n-\alpha)/2).
$$
 Note that $\gamma\to\infty$ as $\alpha\to n$, but $\gamma/(n-\alpha)$ converges to a
 nonzero constant, $\pi^{n/2}2^{n-1}\Gamma(n/2) $. Stein also mentions the Bessel potentials which are associated to the modified inverse Laplace
operators
$$
   ((I -\Delta)^{-\alpha/2}f)  ={\cal I}_{\alpha}(f)
$$
The Bessel potential has a kernel $G_{\alpha}(x)$ that is better
behaved at infinity, though \color{blue} not given by a simple kernel \color{blue} (see \cite{Stein}, page
132).

\medskip

 \subsection*{Fractional Laplacians via  extensions}

In the
case $s=1/2$ it is the well-known technique of harmonic extension to
the upper plane of the space with one more dimension and then taking
the boundary normal derivative, and has been used in the study or
variational problems with thin obstacles as in  \cite{ACld}, see also \cite{Caff1, Caff2}.   For $s\ne 1/2$ the method has been recently developed in \cite{CS}, it involves elliptic equations with weights and weighted normal derivatives.

\section{Comments and extensions}\label{sec-comm}

\noindent {\bf Extensions.}
Very different versions to the evolution process are
obtained when the pressure is related to $u$ in other ways.
Thus, we can use a pressure-density relation of the form $p=f(u)$ with $f$ an increasing function, and then the model equation would be
\begin{equation}
\partial_t u=\nabla\cdot (u\nabla {\cal K}(f(u))).
\end{equation}
Many of the results proved here should apply to this model. \color{blue} Another possibility consists of equations of the form $\partial_t u=\nabla\cdot (f(u)\nabla {\cal K}(u)).$ \normalcolor

\medskip

\noindent {\bf Relaxation. Chemotaxis models.} We could also relax the relation of $p$ to $u$  into the form
\begin{equation}
\partial_t p +(-\Delta)^s p= u.
\end{equation}
This reflection is motivated by a very important system, the Keller-Segel chemotaxis model, \cite{KS, JL}, in which the phenomenon which is modeled by $u$ is not diffusion but
the concentration of a certain population
and $p$ is replaced by variable $c$ proportional to the
concentration of the chemical substance responsible for the
aggregation of the population. A suitable general system is
proposed in the form
\begin{equation}
u_t=\ve \Delta u - \nabla\cdot(u\nabla c), \quad \delta c_t+
{\cal K}^{-1}c=f(u,c).
\end{equation}
The standard chemotaxis model uses ${\cal K}=-\Delta$ and $f(u,c)= u-bc$.
In the limit case where  $\delta$ and $b$ are zero, and if we use
as ${\cal K}$ an integral operator as described above, we get a model
with a term like ours but note the different sign, $u_t=\ve \Delta u -
\nabla\cdot(u\nabla {\cal K}(u))$, which is a consequence of the fact that we
are dealing with aggregation and not diffusion. The study of these equations
is also of interest.

\medskip

\noindent { \bf Finite propagation.} The finite propagation property is not true for other alternative
models of porous medium equation with fractional diffusion like the ones studied in \cite{AC09}, \cite{P4}
The last reference deals with the model
\begin{equation}
\dfrac{\partial u}{\partial t} + (-\Delta)^{1/2} (|u|^{m-1}u)=0,
\end{equation}
For any $m> m_*= (n-1)/n$, it is proved that a unique nonnegative strong solution of this problem exists for data in $L^1_+(\ren)$ and is strictly positive. The maximum principle applies to this problem.

\medskip

\noindent {\bf Uniqueness and comparison.} These are widely open issues. We have used in the present
paper comparison with what we call ``true supersolutions'' and ``true subsolutions''. \color{blue} In one space dimension, a uniqueness proof has been obtained in \cite{BKM} by integrating the equation (with respect to $x$) and using then solutions in the sense of viscosity. Such trick is not available in several dimensions.\normalcolor

\medskip

\noindent {\bf Evolution and regularity of free boundaries.} This is quite important topic motivated
by the property of finite propagation.

\medskip

\noindent {\bf  Smoothing $L^1$ into $L^\infty$ and $C^\alpha$ regularity.} These topics  will be  treated in  \cite{CSV}.

\medskip

\noindent {\bf Asymptotic behaviour.} This topic is under study. Let us outline the main details.
There exists a family of self-similar solutions for this problem, in the spirit of the fundamental solution of the linear problems or the Barenblatt solutions of the standard Porous Medium Equation. The spatial profile of such solutions is obtained by solving for the pressure $p$ an obstacle problem with a truncated paraboloid as obstacle; the corresponding density $u$ is then the mass of the negative fractional Laplacian of $p$, supported on the contact set; all this fits perfectly into the elliptic theory described \cite{ACS}). The details of the  construction, as well as the  convergence of a typical solution to such asymptotic profiles after suitable scaling, will be established in an upcoming publication, \cite{CV2}.

\vskip .5cm


\noindent \textsc{Acknowledgment.} The work was started while JLV  was an Oden Fellow at the ICES
Institute, Univ. of Texas at Austin. LAC was also a guest of the Univ. Aut\'onoma de Madrid. His work was supported by a National Science Foundation grant. JLV was partially supported by Spanish Project MTM2008-06326-C02  and by ESF Programme ``Global and geometric aspects of nonlinear partial differential equations".

\vskip 1cm

\bibliographystyle{amsplain}

\

{\sc Addresses:}

\medskip

{\sc Luis A. Caffarelli}\newline
School of Mathematics, Univ. of Texas at Austin,
1 University Station, C1200, Austin, Texas 78712-1082. \newline
Second affiliation: Institute for Computational Engineering and Sciences.\newline
e-mail: caffarel@math.utexas.edu

\medskip

{\sc Juan Luis V{\'a}zquez}\newline
Departamento de Matem\'{a}ticas, Universidad Aut\'{o}noma de Madrid, 28049
Madrid, Spain.  Second affiliation: Institute ICMAT. \newline
e-mail: juanluis.vazquez@uam.es

\

2000 {\bf Mathematics Subject Classification.} 35K55, 35K65, 76S05.

{\bf Keywords and phases.} Porous medium equation, fractional Laplacian, nonlocal operator, finite propagation.

\end{document}